\documentclass[11pt]{article}
\usepackage{amsmath,amssymb,amscd,enumerate,eucal}
\numberwithin{equation}{section}
\newcommand{\PP}{\mathcal{P}}

\newcommand{\R}{\mathbb{R}}

\newcommand{\Ad}{\text{Ad}}
\newcommand{\Aut}{\text{Aut}}

\newcommand{\YMe}{\mathcal{YM}_{\epsilon}}

\newcommand{\q}{\mathsf{q}}
\newcommand{\A}{\mathsf{a}}

\newcommand{\Exp}{\text{Exp}}
\begin{document}
\setlength{\baselineskip}{17pt}

\title{Small coupling limit and multiple solutions to the Dirichlet
Problem for Yang Mills connections in $4$ dimensions - Part III}
\author{Takeshi Isobe \thanks {%
Tokyo Institute of Technology; email:
isobe.t.ab@m.titech.ac.jp} $\,$ and Antonella Marini \thanks {%
University of L'Aquila / Yeshiva University; email:
marini@dm.univaq.it}}


\date{}
\maketitle

\begin{abstract}
{\footnotesize \baselineskip 4mm} In this paper, the third of this
series, we prove that the spaces $\mathcal{A}^{*,p}_k(A_0;\q)$ and
$\mathcal{B}^{*,p}_{0,k}(A_0;\q)$ which contain $L^p_k$-approximate
solutions to the Dirichlet problem for the $\epsilon$-Yang Mills
equations on a four dimensional disk $B^4$, carry a natural manifold
structure (more precisely a natural structure of Banach bundle), for
$p(k+1)> 4$. All results apply also if $B^4$ is replaced by a
general compact manifold with boundary, and $SU(2)$ is replaced by
any compact Lie group. We also construct bases for the tangent space
to the space of approximate solutions, thus showing that this space
is 8-dimensional for $\epsilon$ sufficiently small, and prove some
technical results used in Parts I and II for the proof of the
existence of multiple solution and, in particular, non-minimal ones,
for this non-compact variational problem.
\end{abstract}

\section{Introduction}
We consider connections $A$ over principal bundles $P$ over the four
dimensional disk, with fiber isomorphic to $SU(2)_\epsilon$, for
$\epsilon >0$, where $SU(2)_\epsilon=SU(2)$ as a set, but it is
endowed with the \emph{$\epsilon$-deformed} multiplication, i.e.,
its Lie algebra is
$(\mathfrak{su}(2)_{\epsilon},[\cdot,\cdot]_{\epsilon})
:=(\mathfrak{su}(2),\epsilon[\cdot,\cdot])$, where $[\cdot,\cdot]$
is the ordinary Lie bracket on $\mathfrak{su}(2)$. We recall that
for a given smooth boundary value $A_0$, the Dirichlet problem for
the $\epsilon$-Yang Mills equations is obtained via a variational
principle from the $SU(2)_{\epsilon}$-Yang Mills functional
\begin{equation}
\label{YMe} \YMe(A)=\int_{B^4}|{F_A}^\epsilon|^2\,dx\,,
\end{equation}
where
${F_A}^\epsilon=dA+\frac{1}{2}[A,A]_{\epsilon}:=dA+\frac{\epsilon}{2}[A,A]$,
and consists of the system $$
(\mathcal{D}_{\epsilon})\quad\qquad\left\{\begin{array}{ll}
{d_A^{\ast}}^\epsilon {F_A}^\epsilon=0\quad&\mbox{in } B^4\\
\iota^{\ast}A\sim A_0\quad&\mbox{at } \partial B^4\,,
\end{array}\right.
$$
where, $\iota:\partial B^4\to\overline{B}^4$ is the inclusion, the
symbol $\sim$  stands for gauge equivalence via a gauge
transformation that extends smoothly to the interior, and
${d_A^{\ast}}^\epsilon:=\ast
d\ast+\ast[A,\ast\cdot]_{\epsilon}:=\ast
d\ast+\epsilon\ast[A,\ast\cdot]$,  where $\ast$ is the Hodge star
operator with respect to the flat metric on $\R^4$.


An absolute minimum, say $\underline{A}_{\epsilon}$, for the
Yang-Mills functional is known to exist by \cite{Marini}. Moreover,
in \cite{IM}, it is shown that the space of connections with
boundary value $A_0$, denoted by $\mathcal{A}(A_0)$, has countable
connected components, i.e.,
$\mathcal{A}(A_0)=\bigsqcup_{j=-\infty}^{\infty}\mathcal{A}_j(A_0)$,
where $\mathcal{A}_j(A_0)$ is the space of connections with relative
2nd Chern number with respect to $\underline{A}_{\epsilon}$ equal to
$j$, and that there always exists a minimizer in
$\mathcal{A}_{+1}(A_0)$ (or $\mathcal{A}_{-1}(A_0)$) if $A_0$ is not
flat. In \cite{IM1, IM2} we have found solutions to $
(\mathcal{D}_{\epsilon})$ in $\mathcal{A}_{+1}(A_0),$ by first
constructing approximate solutions, for small values of the
parameter $\epsilon >0. $ However, some of our proofs utilize
technical results proved in the present paper.

Throughout this paper, we assume $\lambda_0,$ $D_1,$ $D_2,$
$\epsilon$, $\q$ be fixed, with $0<2\lambda_0<d_0$, $0<D_1<D_2$,
$\q:=(p,[g],\lambda)\in\PP(d_0,\lambda_0;D_1,D_2;\epsilon):=\{\q:=(p,[g],\lambda)\in
\PP(d_0,\lambda_0):D_1\epsilon<\lambda^2<D_2\epsilon\},$ where
$\PP(d_0,\lambda_0):= B^4_{1-d_0}\times SO(3)\times(0,\lambda_0)$ is
the parameter space used in \cite{IM1} in the gluing procedure to
construct approximate solutions to the $\epsilon-$Dirichlet problem,
that look like the connected sum of
$\underline{A}_{\epsilon}\#\frac{1}{\epsilon}\text{(1-bubble)}.$

Although all our arguments apply to the general case of $M$, any
smooth compact manifold with boundary, and $G$, any compact Lie
group, we focus on $M=B^4$ and $G= SU(2)$.

We recall (cf. \cite{IM1}) that, for $\q:=(p,g,\lambda)$ fixed,
\begin{equation}
\label{A(q)} A(\q)=\left\{\begin{array}{ll}
(1-\beta_{\lambda,p})\underline{A}_{\epsilon}+\frac{1}{\epsilon}\beta_{\lambda/4,p}\,g\,I^2_{\lambda,p}\,g^{-1}+
\frac{1}{\epsilon}(1-\beta_{\lambda/4,p})gPI^2_{\lambda,p}g^{-1}&\quad\mbox{in } B^4\setminus\{p\}\\
\frac{1}{\epsilon}gI^1_{\lambda,p}g^{-1}&\quad\mbox{in }
B^4_{\lambda/4}(p)\,,
\end{array}\right.
\end{equation}
where
$I^1_{\lambda,p}(x)=\text{Im}\frac{(x-p)d\overline{x}}{\lambda^2+|x-p|^2}$
in $B^4 (p):=\{x\in\R^4:|x-p|<1\}$, together with
$I^2_{\lambda,p}(x)=\text{Im}\frac{\lambda^2(\overline{x}-\overline{p})dx}{|x-p|^2(\lambda^2+|x-p|^2)}$
in $\R^4\setminus\{p\}$ compose the 1-instanton solution to the Yang
Mills equations on $\R^4$,  $\beta(x)=\beta(|x|)\in
C^{\infty}_0(\R^4)$
 is such that $\beta=1$ for $|x|\le 1$, $\beta(x)=0$ for
$|x|\ge 2$ and $0\le\beta(x)\le 1$, and
$\beta_{\lambda,p}(x):=\beta(\lambda^{-1}(x-p)).$

\noindent These connections form the space
\begin{equation}
\label{N}
\mathcal{N}(d_0,\lambda_0):=\{A(\q):\q\in\PP(d_0,\lambda_0)\}\,,
\end{equation}
of approximate solutions to $(\mathcal{D}_{\epsilon})$ with relative
second Chern class equal to +1.

\medskip \noindent Let $g_{12,p}(x)=\frac{x-p}{|x-p|}$ be the
transition maps for the $1$-instanton solution defined above. The
connections $A(\q)$ live on the bundles
\begin{equation}
\label{Pq} P(\q):= \biggl(B^4_{\lambda/4}(p), B^4\setminus\{p\},
g\,g_{12,p}\,g^{-1}\biggr)\,.
\end{equation}
(Notice that $P(p, g, \lambda)= P(p, -g, \lambda)$ and $A(p, g,
\lambda)= A(p, -g, \lambda)$, that's why we take $g\in SU(2)/\{\pm1\}
\cong SO(3)$).

We define $\mathcal{A}^*(A_0;\q)$ as the space of all connections
over $P(\q)$ satisfying $\iota^{\ast}A\sim A_0$ on $\partial B^4$,
via a gauge transformation $g\in \mathcal{G}^*(\q)$ (the group of
all smooth gauge transformations at $\partial B^4$ such that
$g\bigl((1,0,0,0)\bigr)=\mathbf{1}$), which extends smoothly to the
interior, and
$\mathcal{B}(A_0;\q)^*:=\mathcal{A}^*(A_0;\q)/\mathcal{G}^*(\q)$. We
define also the Sobolev counterparts of these spaces:
$\mathcal{A}^{\ast,p}_k(A_0;\q)$, as the space of all connections
over $P(\q)$ satisfying $\iota^{\ast}A\sim A_0$ on $\partial B^4$,
via a gauge transformation $g\in \mathcal{G}^{\ast,p}_{k+1-1/p}(\q)$
(the group of all gauge transformations in $L^p_{k+1-1/p}(\partial
B^4)$ such that $g\bigl((1,0,0,0)\bigr)=\mathbf{1}$), which admits
an $L^p_{k+1}$ extension to the interior, and
$\mathcal{B}^{\ast,p}_k(A_0;\q):=\mathcal{A}^{\ast,p}_k(A_0;\q)/\mathcal{G}^{\ast,p}_{k+1}(\q)$.

\noindent We assume $p(k+1)>4$ throughout the present paper.

\section{A manifold structure for $\mathcal{A}^{*,p}_k(A_0;\q)$ and
$\mathcal{B}^{\ast,p}_k(A_0;\q)$} The space of all connections on a
principal bundle $P$ is an affine space and carries a natural
differentiable structure, while $\mathcal{A}^*(A_0; \q))$ is not an
affine space and it is non-trivial to show that it carries a natural
differentiable structure. The purpose of this section is to show
that $\mathcal{A}^*(A_0; \q)$, its quotient $\mathcal{B}^{\ast}(A_0;
\q)$, and, more in particular, their Sobolev counterparts, possess a
natural differentiable structure.

\noindent We start with  $\mathcal{A}^{\ast,p}_k(A_0; \q)$. In what
follows, we omit the dependence on $\q$ in the notation.
\newtheorem{proposition}{Proposition}[section]
\begin{proposition} \label{P5.1} The space of connections
$\mathcal{A}^{\ast,p}_k(A_0)$  carries a natural structure of Banach
bundle over $\mathcal{G}^{\ast,p}_{k+1-1/p}$. More precisely, for
$A\in\mathcal{A}^{\ast,p}_k(A_0)$ with $\iota^{\ast}A=g^{\ast}A_0$
on $\partial B^4$, for some $g\in\mathcal{G}^{\ast,p}_{k+1-1/p}$,
define
$\pi:\mathcal{A}^{\ast,p}_k(A_0)\to\mathcal{G}^{\ast,p}_{k+1-1/p}$
by $\pi(A)=g$. The map $\pi$ is a well-defined projection, and
$\pi:\mathcal{A}^{\ast,p}_k(A_0)\to\mathcal{G}^{\ast,p}_{k+1-1/p}$
is a vector bundle with fiber isomorphic to
$L^p_{0^T,k}(T^{\ast}B^4\otimes\Ad(P))$.
\end{proposition}
\textit{Proof:} By the Sobolev embedding $L^p_{k+1-1/p}(\partial
B^4)\subset C^0(\partial B^4)$ for $p(k+1)>4$,
$\mathcal{G}^{\ast,p}_{k+1-1/p}$ is a Banach Lie group. Since
$\mathcal{G}^{\ast,p}_{k+1-1/p}$ acts freely on
$\mathcal{A}^{\ast,p}_k(A_0)$ (the restriction of $g$ to the
boundary is uniquely determined), the map
$\pi:\mathcal{A}^{\ast,p}_k(A_0)\to\mathcal{G}^{\ast,p}_{k+1-1/p}$
is well-defined.

\noindent The group $\mathcal{G}^{\ast,p}_{k+1-1/p}$ possesses a
natural manifold structure. In fact, the space
$$L^{\ast,p}_{k+1-1/p}(\Ad(P|_{\partial B^4})):=\{u\in
L^p_{k+1-1/p}(\Ad(P|_{\partial B^4})):u((1,0,0,0))=0\}$$  is a
closed subspace of the Banach space $L^p_{k+1-1/p}(\Aut(P|_{\partial
B^4}))$, which is isomorphic to the Lie algebra of
$\mathcal{G}^{\ast,p}_{k+1-1/p}$. For $u\in
L^{\ast,p}_{k+1-1/p}(\Ad(P|_{\partial B^4}))$, and fixed
$s\in\mathcal{G}^{\ast,p}_{k+1-1/p}$, define $(\Exp_su)(x)=s(x)\exp
u(x)$, where ``$\exp$" denotes the exponential map for the Lie group
$SU(2)$. A local chart of $\mathcal{G}^{\ast,p}_{k+1-1/p}$ at $s$ is
then constructed as follows. Set
\begin{equation}
\label{5.1}
U:=\{u\in L^{\ast,p}_{k+1-1/p}(\Ad(P|_{\partial
B^4})):\|u\|_{k+1-1/p,p}<\sigma\},
\end{equation}
where $0<\sigma<\frac{r}{3}$, where $r$ is the injectivity radius of
$SU(2)$. Let $\tilde{U}_s:=\Exp_s(U)$. With this definition,
$\tilde{U}_s$ is a neighborhood of $s$ in
$\mathcal{G}^{\ast,p}_{k+1-1/p}$ and $\Exp_s:U\to\tilde{U}_s$ yields
a local chart near $s$. See~\cite{FU} for more details.

We construct the `extension operator'
$T:\mathcal{G}^{\ast,p}_{k+1-1/p}\to\mathcal{G}^{\ast,p}_{k+1}$, in
such a way that for any $s\in\mathcal{G}^{\ast,p}_{k+1}$, there
exists a neighborhood $U$  (as in \eqref{5.1}) of $\iota^{\ast}s$ in
$\mathcal{G}^{\ast,p}_{k+1-1/p}$ such that its restriction
$T|U=:T_s:U\to\mathcal{G}^{\ast,p}_{k+1}$ is smooth and satisfies
$T_s(\iota^{\ast}s)=s$.

We take a smooth connection which extends $A_0$ to $P$. For
simplicity, we also denote this connection by $A_0$. Set
$\Delta_{A_0}=\nabla_{A_0}^{\ast}\nabla_{A_0}:L^p_{k+1}(\Ad(P))\to
L^p_{k-1}(\Ad(P))$. Let $U\subset\mathcal{G}^{\ast,p}_{k+1-1/p}$ be
a small neighborhood of $\iota^{\ast}s$. We first define the linear
operator $T_0:U\to L^{\ast,p}_{k+1}(\Ad(P))$ as follows (here we
identify locally $\mathcal{G}^{\ast,p}_{k+1-1/p}$ with its Lie
algebra $L^{\ast,p}_{k+1-1/p}(\Ad(P|_{\partial B^4}))$ via the
exponential map): for given $\xi\in U$, $u=T_0(\xi)$ is the unique
solution of $\Delta_{A_0}u=0$ in $B^4$ with boundary value $u=\xi$
on $\partial B^4$. Elliptic
estimates yield smoothness of $T_0$. Moreover,
$\|T_0(\xi)\|_{k+1,p}\le C\|\xi\|_{k+1-1/p,p}\le C\sigma$, for some
$C>0$ independent of $\xi$. Thus, for small $\sigma$,
$T_s(\xi):=\Exp_s(T_0(\xi))$ is well-defined, smooth with respect to
$s$ and $\xi$, and $T_s(0)=s$.

We next define local charts for $\mathcal{A}^{\ast,p}_k(A_0)$. Let
$U$ be as in \eqref{5.1}. We observe that any
$A\in\mathcal{A}^{\ast,p}_k(A_0)$ can be written as
$A=g^{\ast}(A_0+\alpha)$, for some $g\in\mathcal{G}^{\ast,p}_{k+1}$
and some $\alpha\in L^p_{0^T,k}(T^{\ast}B^4\otimes\Ad(P))$. Hence,
for $s\in\mathcal{G}^{\ast,p}_{k+1-1/p}$, we can define the map
\begin{align}
&\Phi_s:U\times
L^p_{0^T,k}(T^{\ast}\otimes\Ad(P))\to\mathcal{A}^{\ast,p}_k(A_0),\notag\\
&\qquad\Phi_s(u,\alpha)=T_s(u)^{\ast}(A_0+\alpha).\notag
\end{align}

\noindent In order to complete this proof we need the following two
lemmas.

\newtheorem{lemma}{Lemma}[section]
\begin{lemma} \label{L5.1} $\Phi_s$ is one to one.
\end{lemma}
\textit{Proof:} One has that $\Phi_s(u,\alpha)=\Phi_s(v,\beta)$ for
$(u,\alpha),(v,\beta)\in U\times
L^p_{0^T,k}(T^{\ast}B^4\otimes\Ad(P))$ if and only if
$$T_s(u)^{\ast}A_0+T_s(u)^{-1}\alpha T_s(u)=T_s(v)^{\ast}A_0+T_s(v)^{-1}\beta T_s(v)\,.$$
The restriction of this formula to $\partial B^4$ yields
$\iota^{\ast}T_s(u)^{\ast}A_0=\iota^{\ast}T_s(v)^{\ast}A_0$, and,
since the $\mathcal{G}^{\ast,p}_{k+1-1/p}$-action is free,
$\iota^{\ast}T_s(u)=\iota^{\ast}T_s(v)$ on $\partial B^4$. Thus,
$\Exp_su=\Exp_sv$, by definition of $T_s$. From this, $u=v$ and
$\alpha=\beta$. \hfill$\Box$

\medskip

\begin{lemma} \label{L5.2}
$\{\Phi_s,U\times
L^p_{0^T,k}(T^{\ast}B^4\otimes\Ad(P))\}_{s\in\mathcal{G}^{\ast,p}_{k+1-1/p}}$
is a differentiable vector bundle structure for
$\mathcal{A}^{\ast,p}_k(A_0)$.
\end{lemma}
\textit{Proof:} Let $s,t\in\mathcal{G}^{\ast,p}_{k+1-1/p}$ with
$\|s-t\|_{k+1-1/p,p}<\epsilon$ for some small $\epsilon>0$. Here and
in what follows, we replace $\iota^{\ast} s$ by $s$ in our notation.

We need to show that $\Phi_s^{-1}\circ\Phi_t:U\times
L^p_{0^T,k}(T^{\ast}B^4\otimes\Ad(P))\to U\times
L^p_{0^T,k}(T^{\ast}B^4\otimes\Ad(P))$ is smooth for small positive
$\epsilon$. Setting $\Phi_s^{-1}\circ\Phi_t(u,\alpha)=(v,\beta)$,
for $(u,\alpha),(v,\beta)\in U\times
L^p_{0^T,k}(T^{\ast}B^4\otimes\Ad(P))$, one has
$v=\Exp_s^{-1}\circ\Exp_tu$ (cf. Lemma \ref{L5.1}). Then, $\beta$ is given
by
$$\beta=T_s(\Exp_s^{-1}\circ\Exp_tu)\bigl(T_t(u)^{\ast}A_0+T_t(u)^{-1}\alpha
T_t(u)-T_s(\Exp_s^{-1}\circ\Exp_tu)^{\ast}A_0\bigr)\,T_s\,(\Exp_s^{-1}\circ\Exp_tu)^{-1}.$$
Thus, $v$ and $\beta$ depend smoothly on $u$ and $\alpha$. This
completes the proof. \hfill$\Box$

\medskip

\textit{Completion of the proof of Proposition \ref{P5.1}:} We observe that for $s\in\mathcal{G}^{\ast,p}_{k+1-1/p}$
\begin{align*}
\pi^{-1}(s)&=\{A\in\mathcal{A}^{\ast,p}_k(A_0):\iota^{\ast}A=s^{\ast}A_0\}\\
&=\{A={\tilde{s}}^{\ast}(A_0+\alpha):\alpha\in
L^p_{0^T,k}(T^{\ast}B^4\otimes\Ad(P))\}\cong
L^p_{0^T,k}(T^{\ast}B^4\otimes\Ad(P)),
\end{align*}
where $\tilde{s}\in\mathcal{G}^{\ast,p}_{k+1}$ satisfies
$\iota^{\ast}\tilde{s}=s$ and
\begin{align*}
&\pi^{-1}(\tilde{U}_s)=\{A\in\mathcal{A}^{\ast,p}_k(A_0):\iota^{\ast}A=g^{\ast}A_0,~g\in\tilde{U}_s\}\\
&=\{A=T_s(u)^{\ast}(A_0+\alpha):u\in U,~\alpha\in
L^2_{0^T,k}(T^{\ast}B^4\otimes\Ad(P))\} \cong\tilde{U}_s\times
L^p_{0^T,k}(T^{\ast}B^4\otimes\Ad(P)).
\end{align*}
Combining these results with Lemma \ref{L5.2}, one sees that
$\pi:\mathcal{A}^{\ast,p}_k(A_0)\to\mathcal{G}^{\ast,p}_{k+1-1/p}$
is a vector bundle with fiber isomorphic to
$L^p_{0^T,k}(T^{\ast}B^4\otimes\Ad(P))$. This completes the proof.
\hfill$\Box$

\medskip

The following proposition holds for the quotient space
$\mathcal{B}^{\ast,p}_k(A_0)=\mathcal{A}^{\ast,p}_k(A_0)/\mathcal{G}^{\ast,p}_{k+1}$:
\begin{proposition}
\label{P5.2} The space $\mathcal{B}^{\ast,p}_k(A_0)$ defined in $\S1$ has a differentiable manifold structure.
Moreover, the canonical projection
$\mathcal{A}^{\ast,p}_k(A_0)\to\mathcal{B}^{\ast,p}_k(A_0)$ yields a
principal $\mathcal{G}^{\ast,p}_{k+1}$-bundle.
\end{proposition}
\textit{Proof:} The argument is standard (cf.~\cite{DK},~\cite{FU}).
One first seeks a candidate for the slice in
$\mathcal{A}^{\ast,p}_k(A_0)$ at a given
$A_1\in\mathcal{A}^{\ast,p}_k(A_0)$.

For $\xi\in L^{\ast,p}_{k+1}\Ad(P))$ and $t\in\R$ with $|t|$ small,
define $g_t:=\Exp_{\mathbf{1}}(t\xi)\in\mathcal{G}^{\ast,p}_{k+1}$.
There exists $g_0\in\mathcal{G}^{\ast,p}_{k+1}$ and $\alpha_0\in
L^p_{0^T,k}(T^{\ast}B^4\otimes\Ad(P))$ such that
$A_1=g_0^{\ast}(A_0+\alpha_0)$. From
\begin{align*}
g_t^{\ast}A_1&=(\Exp_{g_0}(t\xi))^{\ast}(A_0+\alpha_0)\\
&=
T_{g_0}(t\xi)^{\ast}[A_0+\alpha_0+T_{g_0}(t\xi)(\Exp(t\xi)^{\ast}(A_0+\alpha_0)
-T_{g_0}(t\xi)^{\ast}(A_0+\alpha_0))T_{g_0}(t\xi)^{-1}],
\end{align*}
one has
$$
\Phi_{g_0}^{-1}(g_t^{\ast}A_1)=(t\xi|_{\partial B^4},\alpha_0+
T_{g_0}(t\xi)(\Exp_{g_0}(t\xi)^{\ast}(A_0+\alpha_0)-T_{g_0}(t\xi)^{\ast}(A_0+\alpha_0))T_{g_0}(t\xi)^{-1}).
$$
Thus a tangent vector to the orbit $\{g_t^{\ast}A_1\}_t$ at $A_1$ is
given by
$$\frac{d}{dt}\Big|_{t=0}\Phi_{g_0}^{-1}(g_t^{\ast}A_1)=(\xi|_{\partial B^4},g_0(d_{A_1}\xi-d_{A_1}T_0(\xi))g_0^{-1}).$$
Therefore, a vector $(u,\alpha)\in
L^{\ast,p}_{k+1-1/p}(\Ad(P|_{\partial B^4}))\times
L^p_{0^T,k}(T^{\ast}B^4\otimes\Ad(P))$ in the slice at $A_1$ must
satisfy
\begin{equation}
\label{5.2}
(\xi|_{\partial B^4},u)_{L^2(\partial
B^4)}+(g_0(d_{A_1}(\xi-T_0(\xi))g_0^{-1},\alpha)_{L^2(B^4)}=0\quad\forall
\xi\in L^{\ast,p}_{k+1}(\Ad(P)).
\end{equation}
By taking $\xi|_{\partial B^4}=0$, one has $T_0(\xi)=0$, and \eqref{5.2}
yields
\begin{equation}
\label{5.3}
(g_0(d_{A_1}\xi)g_0^{-1},\alpha)_{L^2(B^4)}=0,
\end{equation}
thus
\begin{equation}
\label{5.4}
d_{A_0+\alpha_0}^{\ast}\alpha=0\quad\text{in $M$}\;.
\end{equation}
On the other hand, \eqref{5.2}, \eqref{5.4} yield $(\xi|_{\partial
B^4},u)_{L^2(\partial B^4)}=0$ for all $\xi\in
L^{\ast,p}_{k+1-1/p}(\Ad(P|_{\partial B^4}))$, thus $u=0$.
Therefore, a suitable candidate for the slice at $A_1$ is
$$
N=\{(0,\alpha)\in L^{\ast,p}_{k+1-1/p}(\Ad(P|_{\partial B^4}))\times
L^p_{0^T,k}(T^{\ast}B^4\otimes\Ad(P)):d_{A_0+\alpha_0}^{\ast}\alpha=0\}.
$$
We must now prove rigorously that $\mathcal{A}^{\ast,p}_k(A_0)$ is
indeed locally diffeomorphic to $N\times\mathcal{G}^{\ast,p}_{k+1}$,
i.e., that  there exists a neighborhood of
$A_1\in\mathcal{A}^{\ast,p}_k(A_0)$, say $\mathcal{O}(A_1)$,  and
$g:\mathcal{O}(A_1)\to\mathcal{G}^{\ast,p}_{k+1}$ such that
$g(A)^{\ast}A\in N$, for all $A\in\mathcal{O}(A_1)$. This way, the
map $\mathcal{A}^{\ast,p}_k(A_0)\ni A\mapsto (g(A)^{\ast},g(A))\in
N\times\mathcal{G}^{\ast,p}_{k+1}$ would provide a local
diffeomorphism. To prove this, we work in local coordinates near
$A_1$, that is, we express the connections
$A\in\mathcal{A}^{\ast,p}_k(A_0)$ near $A_1$ as
$A=T_{g_0}(u)^{\ast}(A_0+\alpha_0+\alpha)$, for some $(u,\alpha)\in
L^{\ast,p}_{k+1-1/p}(\Ad(P|_{\partial B^4}))\times
L^p_{0^T,k}(T^{\ast}B^4\otimes\Ad(P))$, with $\|u\|_{k+1-1/p,p}$
and $\|\alpha\|_{k,p}$ small, and prove the following lemma.
\begin{lemma}
\label{L5.3}
For $(u,\alpha)\in L^{\ast,p}_{k+1-1/p}(\Ad(P|_{\partial
B^4}))\times L^p_{0^T,k}(T^{\ast}B^4\otimes\Ad(P))$, with
$\|u\|_{k+1-1/p,p}$ and $\|\alpha\|_{k,p}$ suitably small, there
exists $g(u,\alpha)\in\mathcal{G}^{\ast,p}_{k+1}$ such that
$g(u,\alpha)^{\ast}A\in N$. Moreover, $g(u,\alpha)$ depends smoothly
on $u$ and $\alpha$, and $g(0,0)=\mathbf{1}$.
\end{lemma}
\textit{Proof:} Since we look for $g$ near the identity, we may
assume that $g$ is of the form $g=\Exp_{\mathbf{1}}\xi$ for $\xi\in
L^{\ast,p}_{k+1}(\Ad(P))$, with $\|\xi\|_{k+1,p}$ small. There
exists $\epsilon>0$ such that the map

\noindent $f:\{u\in L^{\ast,p}_{k+1-1/p}(\Ad(P|_{\partial
B^4})):\|u\|_{k+1-1/p,p}<\epsilon\} \times\{\xi\in
L^{\ast,p}_{k+1}(\Ad(P|_{\partial B^4})):\|\xi\|_{k+1,p}<\epsilon\}
\to L^{\ast,p}_{k+1-1/p}(\Ad(P|_{\partial B^4}))$,
 defined via
$\Exp_{\mathbf{1}}u\cdot\Exp_{\mathbf{1}}\xi=\Exp_{\mathbf{1}}(f(u,\xi))$,
is smooth.

One has
\begin{align*}
g^{\ast}A&=T_{g_0}(f(u,\xi))^{\ast}(A_0+\alpha_0+\alpha)+(T_{g_0}(u)g)^{\ast}(A_0+\alpha_0+\alpha)-T_{g_0}(f(u,\xi))^{\ast}(A_0+\alpha_0+\alpha)\notag\\
&=T_{g_0}(f(u,\xi))^{\ast}(A_0+\alpha_0+\alpha+\beta),
\end{align*}
where
$\beta=\beta(u,\xi)=T_{g_0}(f(u,\xi))[(T_{g_0}(u)g)^{\ast}(A_0+\alpha_0+\alpha)-
T_{g_0}(f(u,\xi))^{\ast}(A_0+\alpha_0+\alpha)]T_{g_0}(f(u,\xi))^{-1}$\;.
Hence,
$$\Phi_{g_0}^{-1}(g^{\ast}A)=(f(u,\xi),\alpha_0+\alpha+\beta(u,\xi)).$$
Thus, we must find $\xi\in L^{\ast,p}_{k+1}(\Ad(P))$ which satisfies
the following equations for small $\|u\|_{k+1-1/p,p}$ and
$\|\alpha\|_{k,p}$:
\begin{align}
&f(u,\xi)=0\quad\text{on $\partial B^4$,}\label{5.5}\\
&d_{A_0+\alpha_0}^{\ast}(\alpha+\beta(u,\xi))=0\quad\text{in $B^4$.}\label{5.6}
\end{align}
The implicit function theorem yields the existence of a solution for
the system \eqref{5.5}-\eqref{5.6} as follows: define
$$F:
L^{\ast,p}_{k+1-1/p}(\Ad(P|_{\partial B^4}))\times
L^p_{0^T,k}(T^{\ast}B^4\otimes\Ad(P))\times
L^{\ast,p}_{k+1}(\Ad(P))\to L^{\ast,p}_{k+1-1/p}(\Ad(P|_{\partial
B^4}))\times L^p_{k-1}(\Ad(P))$$
by
$$F(u,\alpha,\xi):=(f(u,\xi)|_{\partial B^4},d_{A_0+\alpha_0}^{\ast}(\alpha+\beta(u,\xi))),$$
where $g=\Exp_{\mathbf{1}}\xi$. For $\eta\in
L^{\ast,p}_{k+1}(\Ad(P))$, a simple calculation yields
\begin{align*}
\Big\langle\frac{\partial
F}{\partial\xi}(0,0,0),\eta\Big\rangle=\frac{d}{dt}\Big|_{t=0}F(0,0,t\eta)
=(\eta|_{\partial
B^4},d_{A_0+\alpha_0}^{\ast}d_{A_0+\alpha_0}(g_0(\eta-T_0(\eta))g_0^{-1})).
\end{align*}
We claim that $\frac{\partial F}{\partial\xi}(0,0,0):
L^{\ast,p}_{k+1}(\Ad(P))\to L^{\ast,p}_{k+1-1/p}(\Ad(P|_{\partial
B^4}))\times L^p_{k-1}(\Ad(P))$ is an isomorphism. To see this, let
$\eta\in\ker\frac{\partial F}{\partial\xi}(0,0,0)$. Then
$\eta|_{\partial B^4}=0$ and, also,
$d_{A_0+\alpha_0}^{\ast}d_{A_0+\alpha_0}(g_0(\eta-T_0(\eta))g_0^{-1})=0$.
These entail $T_0(\eta)=0$ and $\eta=0$. Thus, $\frac{\partial
F}{\partial\xi}(0,0,0)$ is one-to-one. Furthermore, $\frac{\partial
F}{\partial\xi}(0,0,0)$ is onto, since:
\begin{enumerate}[(i)]
\item For all $\varphi\in L^p_{0,k+1}(\Ad(P))$ and $\eta\in L^{\ast,p}_{k+1}(\Ad(P))$, we have
$T_0((\eta+\varphi)|_{\partial B^4})=T_0(\eta|_{\partial B^4})$.
(Thus, $(\varphi+\eta)-T_0((\varphi+\eta)|_{\partial
B^4})=\varphi+\eta-T_0(\eta|_{\partial B^4})$).
\item $\Ad(g_0^{-1}):L^p_{0,k+1}(\Ad(P))\to L^p_{0,k+1}(\Ad(P))$ is an isomorphism.
\item $\Delta_{A_0+\alpha_0}:=d_{A_0+\alpha_0}^{\ast}d_{A_0+\alpha_0}:L^p_{0,k+1}(\Ad(P))\to L^p_{k-1}(\Ad(P))$ is an isomorphism.
\end{enumerate}
Thus  $\frac{\partial F}{\partial\xi}(0,0,0)$ is an isomorphism and,
 by the implicit function theorem, there exists a neighborhood
$\mathcal{U}$ of $(0,\alpha_0)\in U\times
L^p_{0^T,k}(T^{\ast}\otimes\Ad(P))$ such that for all
$(u,\alpha_0+\alpha)\in\mathcal{U}$, there exists
$g(u,\alpha)\in\mathcal{G}^{\ast,p}_{k+1}$ satisfying
$g(u,\alpha)^{\ast}A\in N$, for
$A=T_{g_0}(u)^{\ast}(A_0+\alpha_0+\alpha)$. This completes the
proof. \hfill$\Box$

\medskip

The next lemma can be proved similarly to the corollary at p. 50 of
\cite{FU}, so we omit its proof.

\begin{lemma} \label{L5.4} Assume $p(k+1)>4$. Then $\mathcal{B}^{\ast,p}_k(A_0)$ is Hausdorff.
\end{lemma}

In order to complete the proof of Proposition \ref{P5.2}, we also need the
following:

\begin{lemma}\label{L5.5}
For any given $(0,\alpha_0)\in N$, there exists $\epsilon>0$ such
that the $\epsilon$-ball of $N$ centered at $(0,\alpha_0)$ injects
into $\mathcal{B}^{\ast,p}_k(A_0)$.
\end{lemma}
\textit{Proof:} Let $A_1=g_0^{\ast}(A_0+\alpha_0)$ (where
$g_0\in\mathcal{G}^{\ast,p}_{k+1}$, and $\alpha_0\in
L^p_{0^T,k}(T^{\ast}B^4\otimes\Ad(P))$). Let
$(0,\alpha_0+\alpha),(0,\alpha_0+\beta)\in N$, with
$\|\alpha\|_{k,p}<\epsilon$ and $\|\beta\|_{k,p}<\epsilon$, where
$\epsilon>0$ is a small number. Set
$A'=g_0^{\ast}(A_0+\alpha_0+\alpha)$,
$A''=g_0^{\ast}(A_0+\alpha_0+\beta)$. Assume that there exists
$g\in\mathcal{G}^{\ast,p}_{k+1}$ such that $g^{\ast}A'=A''$. Under
these assumptions, we need to show that $g=\mathbf{1}$ and
$\alpha=\beta$. To prove this, we rewrite the condition
$g^{\ast}A'=A''$ as
\begin{equation}
\label{5.7}
d_{A_1}g=g(g_0^{-1}\alpha g_0)-(g_0^{-1}\beta g_0)g.
\end{equation}
Let us now write $g:=\tilde{g}+c\in L^p_{k+1}(\text{End}(P))$, where
$c\in\ker d_{A_1}$, $\tilde{g}\in(\ker d_{A_1})^{\perp}$ (here,
$(\ker d_{A_1})^{\perp}$ is the $L^2$-orthogonal complement of $\ker
d_{A_1}$ in $L^p_{k+1}(\text{End}(P)).$ By \eqref{5.7}, there exists
a constant $C>0$ such that
\begin{equation}
\label{5.8}
\|\tilde{g}\|_{p,k+1}\le C\epsilon.
\end{equation}
On the other hand, since $g((0,0,0,1))=\mathbf{1}$,
$c\simeq\mathbf{1}$ and, by \eqref{5.8},
$\|g-\mathbf{1}\|_{k+1,p}<C\epsilon$. By Lemma \ref{L5.3}, one has
$g=\mathbf{1}$ and $\alpha=\beta$. This completes the proof.
\hfill$\Box$

\section{Interaction estimates and construction of orthonormal bases for
$T_{A(\q)}\mathcal{N}(d_0,\lambda_0)$ and
$T_{\tilde{A}(\q)}\tilde{N}(d_0,\lambda_0)$}

In \cite{IM1}, for technical reasons, to the purpose of transforming
the non compact $\epsilon-$Dirichlet problem into a finite
dimensional problem, we introduced the extensions $\tilde{P}(\q)$ to
$\R^4$ of the bundles $P(\q)$,
 \begin{equation}
 \label{Ptilde(q)}
 \tilde{P}(\q):=\biggl(\{B_{\lambda/4}(p), \;\R^4\setminus\{p\}\},
 \;\{gg_{12,p}g^{-1}\}\biggr),\quad \mbox{ for }\,\q\in\PP(d_0,\lambda_0),
 \end{equation}
and the spaces of connections
\begin{equation}
\label{tildeN} \mathcal{\tilde
 N}(d_0,\lambda_0):=\{\tilde{A}(\q):\q\in\PP(d_0,\lambda_0)\}\,,
 \end{equation}
where
\begin{equation}
\label{tildeA(q)} \tilde{A}(\q)=\left\{\begin{array}{ll}
 \frac{1}{\epsilon}gI^2_{\lambda,p}g^{-1}&\quad\mbox{in } \R^4\setminus\{p\}\\
 \frac{1}{\epsilon}gI^1_{\lambda,p}g^{-1}&\quad\mbox{in }
 B^4_{\lambda/4}(p)\,.
 \end{array}\right.
 \end{equation}

 In this section, we prove some technical
lemmas on the tangent space to approximate solutions
$T_{A(\q)}\mathcal{N}(d_0,\lambda_0)$ (cf. $\S1$) and on the tangent
space $T_{\tilde{A}(\q)}\tilde{N}(d_0,\lambda_0).$ These lemmas have
been used in sections $\S3.5$-$\S3.8$ of \cite{IM1}, where, by a
standard technique for non-compact variational problems, the Yang
Mills equation is first solved in
$T_{A(\q)}\mathcal{N}(d_0,\lambda_0)^\perp$ (the orthogonal
complement to the tangent space to approximate solutions),  i.e.,
essentially, orthogonally to the kernel of the Hessian of the
$\epsilon$-Yang Mills functional (cf., in particular, Lemma 3.9
in~\cite{IM1}). This result allowed us to turn the problem into a
finite dimensional one (cf. Proposition 3.2 in \cite{IM1}).

\medskip\noindent
In the following, we choose $\bigl(\xi_1, \xi_2, \xi_3\bigr)$, where
$\xi_1=\begin{pmatrix}
0&0&0\\ 0&0&-1\\ 0&1&0\end{pmatrix}$, $\xi_2=\begin{pmatrix} 0&0&1\\
0&0&0\\ -1&0&0\end{pmatrix}$, $\xi_3=\begin{pmatrix} 0&-1&0\\
1&0&0\\ 0&0&0\end{pmatrix}$, as basis for the Lie algebra
 $\mathfrak{so}(3)$. Right translation by $g$ yields a basis for
$T_{[g]}SO(3)$, which we denote by $\bigl(\xi_1[g], \xi_2[g],
\xi_3[g]\bigr)$.

\begin{lemma} \label{L5.7} The following estimates hold for small $\epsilon>0$:
\begin{align}
&\Big\|\frac{\partial A}{\partial
p_i}(\q)\Big\|_{A(\q);1,2;B^4}\simeq\epsilon^{-3/2}\quad(1\le i\le
4),\\
&\Big\|\frac{\partial
A}{\partial\xi_i[g]}(\q)\Big\|_{A(\q);1,2;B^4}\simeq\epsilon^{-1}\quad(1\le
i\le 3),\\
& \Big\|\frac{\partial
A}{\partial\lambda}(\q)\Big\|_{A(\q);1,2;B^4}\simeq\epsilon^{-3/2},\\
&\Big|\Big(\frac{\partial A}{\partial p_i}(\q),\frac{\partial
A}{\partial
p_j}(\q)\Big)_{A(\q);1,2;B^4}\Big|\lesssim\epsilon^{-3/2}\quad(1\le
i\ne j\le 4),\\
&\Big|\Big(\frac{\partial A}{\partial\xi_i[g]}(\q),\frac{\partial
A}{\partial\xi_j[g]}(\q)\Big)_{A(\q);1,2;B^4}\Big|\lesssim\epsilon^{-1}\quad(1\le
i\ne j\le 3),\\
&\Big|\Big(\frac{\partial A}{\partial p_i}(\q),\frac{\partial
A}{\partial\xi_j[g]}(\q)\Big)_{A(\q);1,2;B^4}\Big|\lesssim\epsilon^{-1}\quad(1\le
i\le 4,~1\le j\le 3),\\
&\Big|\Big(\frac{\partial A}{\partial p_i}(\q),\frac{\partial
A}{\partial\lambda}(\q)\Big)_{A(\q);1,2;B^4}\Big|\lesssim\epsilon^{-2}\quad(1\le
i\le 4),\\
&\Big|\Big(\frac{\partial A}{\partial\xi_i[g]}(\q),\frac{\partial
A}{\partial\lambda}(\q)\Big)_{A(\q);1,2;B^4}\Big|\lesssim\epsilon^{-3/2}\quad(1\le
i\le 3).
\end{align}
\end{lemma}
\textit{Proof:} For  $A(\q)$ in the same gauge as described in $\S1$, one has  (modulo
the action of infinitesimal gauge transformations of $P(\q))$:
\begin{align}
\frac{\partial A}{\partial p_i}(\q)&=-\frac{\partial\beta_{\lambda,p}}{\partial
p_i}\underline{A}_{\epsilon}+
\frac{1}{\epsilon}\frac{\partial\beta_{\lambda/4,p}}{\partial
p_i}gI^2_{\lambda,p}g^{-1}+
\frac{1}{\epsilon}\beta_{\lambda/4,p}g\frac{\partial I^2_{\lambda,p}}{\partial p_i}g^{-1}+\qquad\qquad\qquad\notag\\
&\quad-\frac{1}{\epsilon}\frac{\partial\beta_{\lambda/4,p}}{\partial p_i}gPI^2_{\lambda,p}g^{-1}+
\frac{1}{\epsilon}(1-\beta_{\lambda/4,p})g\frac{\partial PI^2_{\lambda,p}}{\partial p_i}g^{-1}\quad\text{in $B^4\setminus\{p\}$},\label{5.19}\\
\frac{\partial A}{\partial p_i}(\q)&=\frac{1}{\epsilon}g\frac{\partial I^1_{\lambda,p}}{\partial
p_i}g^{-1}\quad\text{in
$B_{\lambda/4}(p)$},\\
\frac{\partial A}{\partial\xi_i[g]}(\q)&=\frac{1}{\epsilon}\beta_{\lambda/4,p}[\xi_i,gI^2_{\lambda,p}g^{-1}]+
\frac{1}{\epsilon}(1-\beta_{\lambda/4,p})[\xi_i,gPI^2_{\lambda,p}g^{-1}]\quad\text{in $B^4\setminus\{p\}$},\\
\frac{\partial A}{\partial\xi_i[g]}(\q)&=\frac{1}{\epsilon}[\xi_i,gI^1_{\lambda,p}g^{-1}]\quad\text{in
$B_{\lambda/4}(p)$},\\
\frac{\partial A}{\partial\lambda}(\q)&=-\frac{\partial\beta_{\lambda,p}}{\partial\lambda}\underline{A}_{\epsilon}+
\frac{1}{\epsilon}\frac{\partial\beta_{\lambda/4,p}}{\partial\lambda}gI^2_{\lambda,p}g^{-1}+
\frac{1}{\epsilon}\beta_{\lambda/4,p}g\frac{\partial I^2_{\lambda,p}}{\partial\lambda}g^{-1}+\qquad\qquad\qquad\notag\\
&\quad-\frac{1}{\epsilon}\frac{\partial\beta_{\lambda/4,p}}{\partial\lambda}gPI^2_{\lambda,p}g^{-1}+
\frac{1}{\epsilon}(1-\beta_{\lambda/4,p})g\frac{\partial PI^2_{\lambda,p}}{\partial\lambda}g^{-1}\quad\text{in $B^4\setminus\{p\}$},\\
\frac{\partial A}{\partial\lambda}(\q)&=\frac{1}{\epsilon}g\frac{\partial
I^1_{\lambda,p}}{\partial\lambda}g^{-1}\quad\text{in
$B_{\lambda/4}(p)$}.
\end{align}
The result follows from these formulas by direct calculation.
\hfill$\Box$

\medskip

\noindent From Lemma \ref{L5.7}, it follows that
$T_{A(\q)}\mathcal{N}(d_0,\lambda_0)$ is an $8$-dimensional linear space, if
$\epsilon>0$ is suitably small.

Let us recall that, for a given connection $A$ on the bundle
$P(\q)$, the inner product $(\alpha, \beta)_{A;1,2;B^4}$ for
one-forms $\alpha$, $\beta$ in $L^2_{1; A}(T^{\ast}
B^4\otimes\Ad(P))$ is defined as:
\begin{equation}
\label{innerprod} (\alpha,
\beta)_{A;1,2;B^4}:=\int_{B^4}\bigl({\nabla_A}^{\epsilon}\alpha,
{\nabla_A}^{\epsilon}\beta\bigr)\,dx + \int_{B^4}\bigl(\alpha,
\beta\bigr)\,dx\
\end{equation}
with ${\nabla_A}^{\epsilon}=\nabla+\epsilon[A,\,\cdot\,]$.

\noindent As a corollary of Lemma \ref{L5.7},
one obtains the following orthonormal basis for $T_{A(\q)}\mathcal{N}(d_0,\lambda_0)$,
the tangent space to approximate solutions.
\begin{lemma}
\label{L5.8}
There exists a basis
$\langle\A_1(\q),\A_2(\q),\ldots,\A_8(\q)\rangle$ for
$T_{A(\q)}\mathcal{N}(d_0,\lambda_0)$, orthonormal with respect to the inner
product $(\cdot,\cdot)_{A(\q);1,2;B^4}$, of the following form:
\begin{align}
\A_1(\q)&=a_{11}(\q)\frac{\partial A}{\partial p_1}(\q),\quad \hbox{
with }\;
a_{11}(\q)\simeq\epsilon^{3/2},\label{5.25}\\
\A_2(\q)&=a_{22}(\q)\frac{\partial A}{\partial
p_2}+a_{21}(\q)\A_1(\q),\quad \hbox{ with }\;
a_{22}(\q)\simeq\epsilon^{3/2},
|a_{21}(\q)|\lesssim\epsilon^{3/2},\\
\A_3(\q)&=a_{33}(\q)\frac{\partial A}{\partial
p_3}(\q)+a_{32}(\q)\A_2(\q)+a_{31}(\q)\A_1(\q),\notag\\
&\quad \hbox{ with }\; a_{33}(\q)\simeq\epsilon^{3/2},~|a_{3i}(\q)|\lesssim\epsilon^{3/2}~(i=1,2),\\
\A_4(\q)&=a_{44}(\q)\frac{\partial A}{\partial
p_4}(\q)+a_{43}(\q)\A_3(\q)+a_{42}(\q)\A_2(\q)+a_{41}(\q)\A_1(\q),\notag\\
&\quad \hbox{ with }\; a_{44}(\q)\simeq\epsilon^{3/2},~
|a_{4i}(\q)|\lesssim\epsilon^{3/2}~(1\le i\le 3),\\
 \A_5(\q)&=a_{55}(\q)\frac{\partial
 A}{\partial\xi_1[g]}(\q)+a_{54}(\q)\A_4(\q)+a_{53}(\q)\A_3(\q)+a_{52}(\q)\A_2(\q)+
 a_{51}(\q)\A_1(\q)\,,\notag\\
 &\quad \hbox{ with }\; a_{55}(\q)\simeq\epsilon,~|a_{5i}(\q)|\lesssim\epsilon^{3/2}~(1\le
 i\le 4),\\
 \A_6(\q)&=a_{66}(\q)\frac{\partial
 A}{\partial\xi_2[g]}(\q)+a_{65}(\q)\A_5(\q)+a_{64}(\q)\A_4(\q)+a_{63}(\q)\A_3(\q)\notag\\
 &\quad+a_{62}(\q)\A_2(\q)+
 a_{61}(\q)\A_1(\q),\notag\\
&\quad \hbox{ with
}a_{66}(\q)\simeq\epsilon,~|a_{65}(\q)|\lesssim\epsilon,~|a_{6i}(\q)|\lesssim\epsilon^{3/2}\quad~(1\le
 i\le 4),\\
 \A_7(\q)&= a_{77}(\q)\frac{\partial
 A}{\partial\xi_3[g]}(\q)+a_{76}(\q)\A_6+a_{75}(\q)\A_(\q)+a_{74}(\q)\A_4(\q)\notag\\
&\quad+ a_{73}(\q)\A_3(\q)+a_{72}(\q)\A_2(\q)+a_{71}(\q)\A_1(\q),\notag\\
&\quad \hbox{ with }\;
a_{77}(\q)\simeq\epsilon,~|a_{76}(\q)|,\,|a_{75}(\q)|\lesssim\epsilon,
 ~|a_{7i}(\q)|\lesssim\epsilon^{3/2}~(1\le i\le 4),\\
 \A_8(\q)&=a_{88}(\q)\frac{\partial
 A}{\partial\lambda}(\q)+a_{87}(\q)\A_7(\q)+a_{86}(\q)\A_6(\q)+a_{85}(\q)\A_5(\q)+a_{84}(\q)\A_4(\q)\notag\\
&\quad+a_{83}(\q)\A_3(\q)+a_{82}(\q)\A_2(\q)+a_{81}(\q)\A_1(\q),\notag\\
&\quad \hbox{ with }\;
a_{88}(\q)\simeq\epsilon^{3/2},~|a_{8i}(\q)|\lesssim\epsilon\quad
 ~(1\le i\le 7).
\end{align}
 \end{lemma}
\textit{Proof:} We apply the Gram-Schmidt's orthogonalization
procedure to construct $\A_i(\q)$ ($1\le i\le 8$) from
$\frac{\partial A}{\partial p_i}(\q)$ ($1\le i\le 4$),
$\frac{\partial A}{\partial\xi_i[g]}(\q)$ ($1\le i\le 3$), and
$\frac{\partial A}{\partial\lambda}(\q)$. The asserted result
follows from Lemma \ref{L5.7} by direct calculation. \hfill$\Box$

\medskip
\noindent We construct the following vector fields $\q_i=\q_i(\q)$ ($1\le i\le
8$) inductively on $\PP(d_0,\lambda_0)$, in such a way that the derivative of $A(\q)$ in the direction $\q_i(\q)$
yields $\A_i(\q)$, i.e., $\A_i(\q)=A_{\q_i}(\q)$ for $1\le i\le
8$. More precisely,
\begin{align}
\q_1(\q)&=a_{11}(\q)\frac{\partial}{\partial p_1},\\
\q_2(\q)&=a_{22}(\q)\frac{\partial}{\partial p_2}+a_{21}\q_1(\q),\\
\q_3(\q)&=a_{33}(\q)\frac{\partial}{\partial
p_3}+a_{32}(\q)\q_2(\q)+a_{31}(\q)\q_1(\q),\\
\q_4(\q)&=a_{44}(\q)\frac{\partial}{\partial
p_4}+a_{43}(\q)\q_3(\q)+a_{42}(\q)\q_2(\q)+a_{41}(\q)\q_1(\q),\\
\q_5(\q)&=a_{55}(\q)\xi_1[g]+a_{54}(\q)\q_4(\q)+a_{53}(\q)\q_3(\q)+a_{52}(\q)\q_2(\q)+a_{51}(\q)\q_1(\q),\\
\q_6(\q)&=a_{66}(\q)\xi_2[g]+a_{65}(\q)\q_5(\q)+a_{64}(\q)\q_4(\q)+a_{63}(\q)\q_3(\q)\notag\\
&\quad+a_{62}(\q)\q_2(\q)+a_{61}(\q)\q_1(\q)\;,\\
\q_7(\q)&=
a_{77}(\q)\xi_3[g]+a_{76}(\q)\q_6(\q)+a_{75}(\q)\q_5(\q)+a_{74}(\q)\q_4+a_{73}(\q)\q_3(\q)\notag\\
&\quad+a_{72}(\q)\q_2(\q)+a_{71}(\q)\q_1(\q)\;,\\
\q_8(\q)&=a_{88}(\q)\frac{\partial}{\partial\lambda}+a_{87}(\q)\q_7(\q)+
a_{86}(\q)\q_6(\q)+a_{85}(\q)\q_5(\q)+a_{84}(\q)\q_4(\q)\notag\\
&\quad+a_{83}(\q)\q_3+a_{82}(\q)\q_2(\q)+a_{81}(\q)\q_1(\q).
\end{align}

\noindent Let now $\tilde{\A}_i(\q)=\tilde{A}_{\q_i}(\q)$ (for $1\le
i\le 8$), where $\tilde{A}(\q)$ are the instanton solutions
introduced at the beginning of this section. In the following lemma,
we construct a special basis for
$T_{\tilde{A}(\q)}\tilde{N}(d_0,\lambda_0),$ orthonormal with
respect to the weighted inner product on
$L^2_{1;\tilde{A}(\q)}(T^{\ast}\R^4\otimes\Ad(\tilde{P}))$ defined
by
\begin{equation}
\label{wip}
(\alpha,\beta)_{1,2;\tilde{A}(\q)}:=\int_{\R^4}\bigl({\nabla_{\tilde{A}(\q)}}^{\epsilon}\alpha,{\nabla_{\tilde{A}(\q)}}^{\epsilon}\beta\bigr)\,dx
+ \int_{\R^4} w(x)\bigl(\alpha,\beta\bigr)\,dx,
\end{equation}
where $w(x)=1$ for $|x|\le 1$, and $w(x)=1/(1+|x|^2)^2$ for $|x|>
1$.

\begin{lemma}
\label{L5.9} There exists a basis $\langle\hat{\A}_1(\q),\hat{\A}_2(\q),\ldots,\hat{\A}_8(\q)\rangle$ for
$T_{\tilde{A}(\q)}\tilde{N}(d_0,\lambda_0)$, orthonormal with
respect to the inner product
$(\cdot,\cdot)_{\tilde{A}(\q);1,2;\R^4}$, of the following form:
\begin{align}
\hat{\A}_1(\q)&=b_{11}(\q)\tilde{\A}_1(\q),\quad \hbox{ with }\;
|b_{11}(\q)-1|\lesssim\epsilon,\\
\hat{\A}_2(\q)&=b_{22}(\q)\tilde{\A}_2(\q)+
b_{21}(\q)\hat{\A}_1(\q),\quad \hbox{ with }\;|b_{22}(\q)-1|\lesssim\epsilon,~|b_{21}(\q)|\lesssim\epsilon,\\
\hat{\A}_3(\q)&=b_{33}(\q)\tilde{\A}_3(\q)+b_{32}(\q)\hat{\A}_2(\q)+b_{31}(\q)\hat{\A}_1(\q),\notag\\
&\quad \hbox{ with }\;|b_{33}(\q)-1|\lesssim\epsilon,~|b_{3i}(\q)|\lesssim\epsilon~(i=1,2),\\
\hat{\A}_4(\q)&=b_{44}(\q)\tilde{\A}_4(\q)+b_{43}(\q)\hat{\A}_3(\q)+b_{42}(\q)\hat{\A}_2(\q)+
b_{41}(\q)\hat{\A}_1(\q),\notag\\
&\quad \hbox{ with
}|b_{44}(\q)-1|\lesssim\epsilon,~|b_{4i}(\q)|\lesssim\epsilon~(1\le
i\le 3),\\
\hat{\A}_5(\q)&=b_{55}(\q)\tilde{\A}_5(\q)+b_{54}(\q)\hat{\A}_4(\q)+b_{53}(\q)\hat{\A}_3(\q)+
b_{52}(\q)\hat{\A}_2(\q)+b_{51}(\q)\hat{\A}_1(\q),\notag\\
&\quad \hbox{ with
}|b_{55}(\q)-1|\lesssim\epsilon,~|b_{54}(\q)|\lesssim\epsilon~(1\le
i\le 4),\\
\hat{\A}_6(\q)&=b_{66}(\q)\tilde{\A}_6(\q)+b_{65}(\q)\hat{\A}_5(\q)+b_{64}(\q)\hat{\A}_4(\q)+
b_{63}(\q)\hat{\A}_3(\q)\notag\\
&\quad+b_{62}(\q)\hat{A}_2(\q)+ b_{61}(\q)\hat{\A}(\q),~ \hbox {
with
}\;|b_{66}(\q)-1|\lesssim\epsilon,~|b_{6i}(\q)|\lesssim\epsilon~(1\le
i\le 5),\\
\hat{\A}_7(\q)&=b_{77}(\q)\tilde{\A}_7(\q)+b_{76}(\q)\hat{\A}_6(\q)+b_{75}(\q)\hat{\A}_5(\q)+b_{74}(\q)\hat{\A}_4(\q)+b_{73}(\q)\hat{\A}_3(\q)
\notag\\
&\quad+b_{72}(\q)\hat{\A}_2(\q) +b_{71}(\q)\hat{\A}(\q),~\hbox {
with
}\;|b_{77}(\q)-1|\lesssim\epsilon,~|b_{7i}(\q)|\lesssim\epsilon~(1\le
i\le 6),\\
\hat{\A}_8(\q)&=b_{88}(\q)\tilde{\A}_8(\q)+b_{87}(\q)\hat{\A}_7(\q)+b_{86}(\q)\hat{\A}_6(\q)+b_{85}(\q)\hat{\A}_5(\q)+
b_{84}(\q)\hat{\A}_4(\q)\notag\\
&\quad+b_{83}(\q)\hat{\A}_3(\q)
+b_{82}(\q)\hat{\A}_2(\q)+b_{81}(\q)\hat{\A}_1(\q),~\hbox { with
}\;|b_{88}(\q)-1|\lesssim\epsilon,~|b_{8i}(\q)|\lesssim\epsilon~(1\le
i\le 7).
\end{align}
\end{lemma}
\textit{Proof:} We apply the Gram-Schmidt's orthogonalization
procedure to construct $\hat{\A}_i(\q)$ ($1\le i\le 8$) from
$\tilde{\A}_i(\q)$ ($1\le i\le 8$). In the gauge used for $\tilde
A(\q)$ (cf. beginning of $\S3$), the derivatives of $\tilde{A}(\q)$
are written as (modulo the action of the infinitesimal gauge
transformation of $\tilde{P}(\q)$):
\begin{align}
\frac{\partial\tilde{A}}{\partial p_i}(\q)=
\frac{1}{\epsilon}g\frac{\partial I^2_{\lambda,p}}{\partial
p_i}g^{-1}\quad\text{in $\R^4\setminus\{p\}$},
\quad\frac{\partial\tilde{A}}{\partial
p_i}(\q)=\frac{1}{\epsilon}g\frac{\partial I^1_{\lambda,p}}{\partial
p_i}g^{-1} \quad\text{in $B_{\lambda/4}(p)$},\label{5.49}\\
\frac{\partial\tilde{A}}{\partial\xi_i[g]}(\q)=\frac{1}{\epsilon}[\xi_i,gI^2_{\lambda,p}g^{-1}]
\quad\text{in
$\R^4\setminus\{p\}$},\quad\frac{\partial\tilde{A}}{\partial\xi_i[g]}(\q)=\frac{1}{\epsilon}[\xi_i,gI^1_{\lambda,p}g^{-1}]
\quad\text{in $B_{\lambda/4}(p)$},\label{5.50}\\
\frac{\partial\tilde{A}}{\partial\lambda}(\q)=\frac{1}{\epsilon}\frac{\partial\tilde{A}}{\partial\lambda}(\q)
\quad\text{in
$\R^4\setminus\{p\}$},\quad\frac{\partial\tilde{A}}{\partial\lambda}(\q)=
\frac{1}{\epsilon}\frac{\partial\tilde{A}}{\partial\lambda}(\q)\quad\text{in
$B_{\lambda/4}(p)$}\label{5.51}
\end{align}
From \eqref{5.49}-\eqref{5.51}, by direct computation, one has
\begin{equation}
\label{5.52}
\|\tilde{\A}_i(\q)\|_{\tilde{A}(\q);1,2;\R^4\setminus B^4}\lesssim\epsilon^{3/2}~(1\le i\le 4),\quad
\|\tilde{\A}_i(\q)\|_{\tilde{A}(\q);1,2;\R^4\setminus B^4}\lesssim\epsilon~(5\le i\le 8).
\end{equation}
On the other hand,
\begin{equation}
\label{5.53}
(\tilde{\A}_i(\q),\tilde{\A}_j(\q))_{\tilde{A}(\q);1,2}=\delta_{ij}+O(\epsilon),
\end{equation}
for $1\le i,j\le 8$. In fact,
\begin{align}
\label{5.54}
(\tilde{\A}_i(\q),\tilde{\A}_j(\q))_{\tilde{A}(\q);1,2}&=(\tilde{\A}_i(\q),\tilde{\A}_j(\q))_{\tilde{A}(\q);1,2;B^4}+
(\tilde{\A}_i(\q),\tilde{\A}_j(\q))_{\tilde{A}(\q);1,2;\R^4\setminus B^4}\notag\\
&=(\tilde{\A}_i(\q),\tilde{\A}_j(\q))_{\tilde{A}(\q);1,2;B^4}+O(\epsilon^2)~(\text{by (3.47)})\notag\\
&=(\A_i(\q)+\tilde{\A}_i(\q)-\A_i(\q),\A_j(\q)+\tilde{\A}_j(\q)-\A_j(\q))_{\tilde{A}(\q);1,2;B^4}+O(\epsilon^2)\notag\\
&=(\A_i(\q),\A_j(\q))_{\tilde{A}(\q);1,2;B^4}+O(\epsilon)~(\text{by
Lemma \ref{L3.6} below}).
\end{align}
Here, by setting $b(\q):=\tilde{A}(\q)-A(\q)$, one obtains
\begin{align}
\label{5.55}
(\A_i(\q),\A_j(\q))_{\tilde{A}(\q);1,2;B^4}&=(\A_i(\q),\A_j(\q))_{A(\q);1,2;B^4}+
\int_{B^4}({\nabla_{\tilde{A}(\q)}}^{\epsilon}\A_i(\q),{\nabla_{\tilde{A}(\q)}}^{\epsilon}\A_j(\q))\,dx\notag\\
&\qquad-\int_{B^4}({\nabla_{A(\q)}}^{\epsilon}\A_i(\q),{\nabla_{A(\q)}}^{\epsilon}\A_j(\q))\,dx\notag\\
&=\delta_{ij}+\epsilon\int_{B^4}({\nabla_{A(\q)}}^{\epsilon}\A_i(\q),[b(\q),\A_j(\q)])\,dx+
\epsilon\int_{B^4}([b(\q),\A_i(\q)],{\nabla_{A(\q)}}^{\epsilon}\A_j(\q))\,dx\notag\\
&\qquad+\epsilon^2\int_{B^4}([b(\q),\A_i(\q)],[b(\q),\A_j(\q)])\,dx
=\delta_{ij}+O(\epsilon),
\end{align}
since $\epsilon\|b(\q)\|_{\infty}\lesssim\epsilon$.

Combining \eqref{5.54}, \eqref{5.55}, one obtains \eqref{5.53}, hence the
$\hat{\A}_i(\q)$'s ($1\le i\le 8$) satisfy the assertion of the
lemma. \hfill$\Box$

\medskip
\begin{lemma}
\label{L3.6}
For $\q\in\PP(d_0,\lambda_0;D_1,D_2;\epsilon)$,
\begin{align*}\|\A_i(\q)-\tilde{\A}_i(\q)\|_{A(\q);1,2;B^4}\lesssim\epsilon^{3/2}\quad(1\le
i\le 4),\\
\|\A_i(\q)-\tilde{\A}_i(\q)\|_{A(\q);1,2;B^4}\lesssim\epsilon\quad(5\le
i\le 8).\end{align*}
\end{lemma}

\textit{Proof:} The connections $A(\q)$ and $\tilde{A}(\q)$,
represented in the same gauge as in the previous lemmas,
 differ only on the domain $B^4\setminus
B_{\lambda/4}(p)$. In particular,
$A(\q)-\tilde{A}(\q)=(\beta_{\lambda,p})\underline{A}_{\epsilon}+
\frac{1}{\epsilon}(\beta_{\lambda/4,p}-1)gh_{\lambda,p}g^{-1}$. One
has (modulo infinitesimal gauge transformations):
\begin{align}
\frac{\partial A}{\partial
p_i}(\q)-\frac{\partial\tilde{A}}{\partial p_i}(\q)&=-
\frac{\partial\beta_{\lambda,p}}{\partial
p_i}\underline{A}_{\epsilon}+
\frac{1}{\epsilon}\frac{\partial\beta_{\lambda/4,p}}{\partial
p_i}gh_{\lambda,p}g^{-1}+
\frac{1}{\epsilon}(\beta_{\lambda/4,p}-1)g\frac{\partial
h_{\lambda,p}}{\partial p_i}g^{-1},\\
\frac{\partial
A}{\partial\xi_i[g]}(\q)-\frac{\partial\tilde{A}}{\partial\xi_i[g]}(\q)&=
\frac{1}{\epsilon}(\beta_{\lambda/4,p}-1)[\xi_i,gh_{\lambda,p}g^{-1}],\\
\frac{\partial
A}{\partial\lambda}(\q)-\frac{\partial\tilde{A}}{\partial\lambda}(\q)&=
-\frac{\partial\beta_{\lambda,p}}{\partial\lambda}\underline{A}_{\epsilon}+\frac{1}{\epsilon}\frac{\partial\beta_{\lambda/4,p}}{\partial\lambda}gh_{\lambda,p}g^{-1}+\frac{1}{\epsilon}(\beta_{\lambda/4,p}-1)g\frac{\partial
h_{\lambda,p}}{\partial\lambda}g^{-1}.
\end{align}
The asserted result follows from Lemma \ref{L5.8} by direct computation.
\hfill$\Box$

\medskip
We recall that the Hessian of $\YMe$, denoted by $\nabla^2\YMe(A)$,
is defined by
\begin{equation}
\label{3.58}
\langle\nabla^2\YMe(A)a,b\rangle:=2\int_{B^4}({d_A}^{\epsilon}a,{d_A}^{\epsilon}b)
+2\int_{B^4}(F_A,\epsilon[a,b])\;\hbox{ for all } a, b\in
L^2_{1,0}(T^{\ast}B^4\otimes\Ad(P))\,,
\end{equation}
 where
$\langle\cdot,\cdot\rangle$ denotes the pairing between
$L^2_{1,0}(T^{\ast}B^4\otimes\Ad(P))$ and its dual and ${d_A}^{\epsilon}=d+\epsilon[A,\,\cdot\,]$.
\begin{lemma}
\label{L3.7}
For $\q\in\PP(d_0,\lambda_0;D_1,D_2;\epsilon)$,
\begin{align*}\|\big(\nabla^2\YMe(A(\q))-\nabla^2\YMe(\tilde{A}(\q))\bigr)\A_i(\q)\|_{A(\q);1,2,\ast}&
\lesssim\epsilon^{3/2}\quad(1\le
i\le 8),\\
\|(d_{A(\q)}d_{A(\q)}^{\ast}-d_{\tilde{A}(\q)}d_{\tilde{A}(\q)}^{\ast})\A_i(\q)\|_{A(\q);1,2,\ast}&
\lesssim\epsilon^{3/2}\quad(1\le
i\le 8).
\end{align*}
\end{lemma}

\textit{Proof:} We perform the calculation for $i=1$ (the
remaining cases are analogous).  For $b(\q):=\tilde{A}(\q)-A(\q)$,
$\alpha\in L^2_1(T^{\ast}B^4\otimes\Ad(P(\q)))$, $\beta\in
L^2_{0,1}(T^{\ast}B^4\otimes\Ad(P(\q)))$, one has
\begin{align}
\label{5.59}
&\frac{1}{2}(\nabla^2\YMe(\tilde{A}(\q))-\nabla^2\YMe(A(\q))(\alpha,\beta)\notag\\
&=\epsilon\int_{B^4}({d_{A(\q)}}^{\epsilon}\alpha,[b(\q),\beta])\,dx+\epsilon\int_{B^4}([b(\q),\alpha],{d_{A(\q)}}^{\epsilon}\beta)\,dx\notag\\
&\quad+\epsilon\int_{B^4}({d_{A(\q)}}^{\epsilon}b(\q),[\alpha,\beta])\,dx+\epsilon^2\int_{B^4}([b(\q),\alpha],[b(\q),\beta])\,dx\notag\\
&\quad+\frac{\epsilon^2}{2}\int_{B^4}([b(\q),b(\q)],[\alpha,\beta])\,dx.
\end{align}
We set $\alpha=\A_1(\q)$ in \eqref{5.59}, and estimate each term. In the
following, we write $r:=|x-p|$. Since
\begin{equation}
\label{5.60}
|\A_1(\q)|=a_{11}(\q)\frac{\partial A}{\partial p_i}(\q)
\lesssim\epsilon|\nabla\beta(\lambda^{-1}(\cdot-p))|+\frac{\lambda^3}{r^2(\lambda^2+r^2)}+|\nabla\beta(\lambda^{-1}(\cdot-p))|\lambda^2
+\lambda^3
\end{equation}
and
\begin{equation}
\label{5.61}
|\epsilon
{d_{A(\q)}^{\ast}}^{\epsilon}[b(\q),\beta]|\lesssim(\epsilon\lambda^{-1}|\nabla(\lambda^{-1}(\cdot-p))|+\lambda^2+\epsilon|I^2_{\lambda,p}|)|\beta|+
\epsilon|\nabla\beta|,
\end{equation}
one obtains
\begin{align}
\label{5.62}
&\biggl|\epsilon\int_{B^4}({d_{A(\q)}}^{\epsilon}\alpha,[b(\q),\beta]) \biggr|=\biggl|\epsilon\int_{B^4}(\alpha,{d_{A(\q)}^{\ast}}^{\epsilon}[b(\q),\beta]) \biggr|\notag\\
&\lesssim\int_{B_{2\lambda}(p)\setminus
B_{\lambda}(p)}\Big(\epsilon^2\lambda^{-1}+\epsilon\lambda+\lambda^4+
\frac{\epsilon^2\lambda^2}{r(\lambda^2+r^2)}+\frac{\epsilon\lambda^2}{r^2(\lambda^2+r^2)}+
\frac{\epsilon\lambda^4}{r(\lambda^2+r^2)}\Big)|\beta|\,dx\notag\\
&\quad+\int_{B_{2\lambda}(p)\setminus B_{\lambda}(p)}(\epsilon^2+\epsilon\lambda^2)|\nabla\beta|\,dx\notag\\
&\quad+\int_{B^4\setminus B_{\lambda/4}(p)}\Big(\lambda^5+\frac{\lambda^5}{r^2(\lambda^2+r^2)}+\frac{\epsilon\lambda^5}{r^3(\lambda^2+r^2)^2}+
\frac{\epsilon\lambda^5}{r(\lambda^2+r^2)}\Big)|\beta|\,dx\notag\\
&\quad+\int_{B^4\setminus B_{\lambda/4}(p)}\Big(\epsilon\lambda^3+\frac{\epsilon\lambda^3}{r^2(\lambda^2+r^2)}\Big)|\nabla\beta|\,dx\notag\\
&\lesssim (\epsilon^2\lambda^2+\epsilon\lambda^4+\lambda^7+\epsilon^2\lambda^2+\epsilon\lambda+\epsilon\lambda^4)
\|\beta\|_{L^4(B_{2\lambda}(p)\setminus B_{\lambda}(p))}+(\epsilon^2\lambda^2+\epsilon\lambda^4)\|\nabla\beta\|_{L^2(B_{2\lambda}(p)\setminus B_{\lambda}(p))}\notag\\
&\quad+(\lambda^5+\lambda^4+\epsilon\lambda+\epsilon\lambda^5|\log\lambda|)\|\beta\|_{L^4(B^4\setminus B_{\lambda/4}(p))}+(\epsilon\lambda^3+\epsilon\lambda^2)\|\nabla\beta\|_{L^2(B^4\setminus B_{\lambda/4}(p))}\notag\\
&\lesssim\epsilon^{3/2}\|\beta\|_{A(\q);1,2}.\qquad\qquad\qquad
\end{align}
From \eqref{5.62}, the first part of Lemma \ref{L3.7} (for $i=1$) follows.

To prove the second part, for any given $\beta\in
L^2_{0,1}(T^{\ast}B^4\otimes\Ad(P(\q)))$ we estimate
\begin{align}
\label{5.63}
&\langle({d_{A(\q)}}^{\epsilon}{d_{A(\q)}^{\ast}}^{\epsilon}-{d_{\tilde{A}(\q)}}^{\epsilon}{d_{\tilde{A}(\q)}^{\ast}}^{\epsilon})\A_i(\q),\beta\rangle
=\int_{B^4}({d_{A(\q)}^{\ast}}^{\epsilon}\A_i(\q),{d_{A(\q)}^{\ast}}^{\epsilon}\beta)\,dx-
\int_{B^4}({d_{\tilde{A}(\q)}^{\ast}}^{\epsilon}\A_i(\q),{d_{\tilde{A}(\q)}^{\ast}}^{\epsilon}\beta)\,dx\notag\\
&=\int_{B^4}(({d_{A(\q)}^{\ast}}^{\epsilon}-{d_{\tilde{A}(\q)}^{\ast}}^{\epsilon})\A_i(\q),{d_{A(\q)}^{\ast}}^{\epsilon}\beta)\,dx+
\int_{B^4}({d_{\tilde{A}(\q)}^{\ast}}^{\epsilon}\A_i(\q),({d_{A(\q)}^{\ast}}^{\epsilon}-{d_{\tilde{A}(\q)}^{\ast}}^{\epsilon})\beta)\,dx\notag\\
&=\int_{B^4\setminus B_{\lambda/4}(p)}(\epsilon\ast[b(\q),\ast\A_i(\q)],{d_{A(\q)}^{\ast}}^{\epsilon}\beta)\,dx+
\int_{B^4\setminus B_{\lambda/4}(p)}({d_{\tilde{A}(\q)}^{\ast}}^{\epsilon}\A_i(\q),\epsilon\ast[b(\q),\ast\beta])\,dx\notag\\
&\lesssim\epsilon\biggl(\int_{B^4\setminus
B_{\lambda/4}(p)}|\A_i(\q)|^2\,dx\biggr)^{1/2}\|\beta\|_{A(\q);1,2}+\epsilon\biggl(\int_{B^4\setminus
B_{\lambda/4}(p)}|{d_{\tilde{A}(\q)}^{\ast}}^{\epsilon}\A_i(\q)|^{4/3}\,dx\biggr)^{3/4}\|\beta\|_{A(\q);1,2}.
\end{align}
By \eqref{5.60}, the first integral in \eqref{5.63} is estimated as
\begin{align}
\label{5.64}
\int_{B^4\setminus B_{\lambda/4}(p)}|\A_1(\q)|^2\,dx&\lesssim\epsilon^2\lambda^4+\lambda^2+\lambda^8+\lambda^6\lesssim\epsilon\,.
\end{align}
As for the second integral, since
$|{d_{\tilde{A}(\q)}^{\ast}}^{\epsilon}\A_1|\le|d^{\ast}\A_1|+\epsilon|[\tilde{A}(\q),\ast\A_1(\q)]|$,
from \eqref{5.19}, \eqref{5.25}, \eqref{5.49}, \eqref{5.60}, it follows that
\begin{align}
\label{5.65}
|{d_{\tilde{A}(\q)}^{\ast}}^{\epsilon}\A_1(\q)|&\lesssim\epsilon\lambda^{-1}|\nabla^2\beta(\lambda(\cdot-p))|+
\epsilon|\nabla\beta(\lambda(\cdot-p))|+
\frac{\lambda^3}{r^3(\lambda^3+r^2)}\notag\\
&\quad+\lambda|\nabla^2\beta(\lambda/4(\cdot-p))|+\lambda^2|\nabla\beta(\lambda(\cdot-p))|+\lambda^3\notag\\
&\quad+\frac{\lambda^2}{r(\lambda^2+r^2)}\Big(\epsilon|\nabla\beta(\lambda(\cdot-p))|+
\frac{\lambda^3}{r^2(\lambda^2+r^2)}+
\lambda^2|\nabla\beta(\lambda/4(\cdot-p))|+\lambda^3\Big).
\end{align}
Thus,
\begin{align}
\label{5.66}
&\int_{B^4\setminus B_{\lambda/4}(p)}|{d_{\tilde{A}(\q)}^{\ast}}^{\epsilon}\A_1(\q)|^{4/3}\,dx\lesssim\epsilon^{4/3}\lambda^{8/3}+
\epsilon^{4/3}\lambda^4+\int_{\lambda/4}^2\frac{\lambda^4}{r(\lambda^2+r^2)^{4/3}}\,dr\notag\\
&\quad+\lambda^{16/3}+\lambda^{20/3}+\lambda^4+
\epsilon^{4/3}\lambda^{8/3}\int_{\lambda}^{2\lambda}\frac{r^{5/3}}{(\lambda^2+r^2)^{4/3}}\,dr
+
\lambda^{20/3}\int_{\lambda/4}^2\frac{1}{r(\lambda^2+r^2)^{8/3}}\,dr\notag\\
&\quad+\lambda^{16/3}\int_{\lambda}^{2\lambda}\frac{r^{5/3}}{(\lambda^2+r^2)^{4/3}}\,dr+
\lambda^{20/3}\int_{\lambda/4}^2\frac{r^{5/3}}{(\lambda^2+r^2)^{4/3}}\,dr\notag\\
&\lesssim\,\epsilon^{4/3}\lambda^{8/3}+\epsilon^{4/3}\lambda^4+\lambda^{4/3}+\lambda^{16/3}+\lambda^{20/3}+
\lambda^4
+\epsilon^{4/3}\lambda^{8/3}+\lambda^{4/3}+\lambda^{16/3}+\lambda^{20/3}|\log\lambda|
\lesssim \,\epsilon^{2/3}.
\end{align}
Combining \eqref{5.63}, \eqref{5.64}, \eqref{5.66}, one obtains
\begin{equation}
\label{5.67}
|\langle({d_{A(\q)}}^{\epsilon}{d_{A(\q)}^{\ast}}^{\epsilon}-{d_{\tilde{A}(\q)}}^{\epsilon}{d_{\tilde{A}(\q)}^{\ast}}^{\epsilon})\A_i(\q),\beta\rangle|
\lesssim\epsilon^{3/2}\|\beta\|_{A(\q);1,2},
\end{equation}
i.e., the second assertion of Lemma \ref{L3.7} for $i=1$. \hfill$\Box$

\medskip
\begin{lemma} \label{L3.10} Let $\A_i(\q)$ be the elements of the orthonormal basis constructed in Lemma \ref{L5.8}.  The following estimates hold:
\begin{equation}
\label{140}
\|{\A_i}_{\q_j}(\q)^{\perp}\|_{A(\q);1,2;B^4}\lesssim\epsilon
\end{equation}
for $1\le i,j\le 8$, where ${\A_i}_{\q_j}(\q)$ denotes the
directional derivative of $\A_i(\q)$ in the direction $\q_j$.
\end{lemma}
\textit{Proof:} We prove the lemma for $i=1$, $j=1$.
The remaining cases are similar. For $A(\q)$ as represented  in $\S1$, we first observe that
$$({\A_1}_{\q_1}(\q))^{\perp}=
\Big(a_{11}(\q)\frac{\partial}{\partial p_1}\Big(a_{11}(\q)\frac{\partial A}{\partial p_1}(\q)\Big)\Big)^{\perp}
=a_{11}(\q)^2\Big(\frac{\partial^2 A}{\partial
p_1^2}(\q)\Big)^{\perp}\,$$
so, we need to estimate
$\Big\|a_{11}(\q)^2\frac{\partial^2 A}{\partial
p_1^2}(\q)\Big\|_{A(\q);1,2}$.

On $B^4\setminus B_{\lambda/4}(p)$,
\begin{align}
\label{5.68}
\epsilon\frac{\partial^2 A}{\partial p_1^2}(\q)&=
\,-\epsilon\lambda^{-2}\frac{\partial^2\beta}{\partial x_1^2}(\lambda^{-1}(\cdot-p))\underline{A}_{\epsilon}+
g\frac{\partial^2 I^2_{\lambda,p}}{\partial p_1^2}g^{-1}+
16\lambda^{-2}\frac{\partial^2\beta}{\partial x_1^2}(4\lambda^{-1}(\cdot-p))gh_{\lambda,p}g^{-1}\notag\\
&\quad-4\lambda^{-1}\frac{\partial\beta}{\partial
x_1}(4\lambda^{-1}(\cdot-p))g\frac{\partial h_{\lambda,p}}{\partial
p_1}g^{-1}+(\beta_{\lambda/4,p}-1)g\frac{\partial^2
h_{\lambda,p}}{\partial p_1^2}g^{-1}.
\end{align}
Hence,
\begin{align}
\label{5.69}
\Big|\epsilon\frac{\partial^2 A}{\partial p_1^2}(\q)&\Big|\lesssim\epsilon\lambda^{-2}|\nabla^2\beta(\lambda^{-1}(\cdot-p))|+
\frac{\lambda^2}{r^3(\lambda^2+r^2)}+
|\nabla^2\beta(\lambda^{-1}(\cdot-p))|\notag\\
&\quad+\lambda|\nabla\beta(\lambda^{-1}(\cdot-p))|+\lambda^2,
\end{align}
\begin{align}
\label{5.70}
\Big|\epsilon\nabla\Big(\frac{\partial^2 A}{\partial p_1^2}(\q)\Big)&\Big|\lesssim\epsilon\lambda^{-3}|\nabla^3\beta(\lambda^{-1}(\cdot-p))|+
\epsilon\lambda^{-2}|\nabla^2\beta(\lambda^{-1}(\cdot-p))|+\frac{\lambda^2}{r^4(\lambda^2+r^2)}\notag\\
&\quad+\lambda^{-1}|\nabla^3\beta(\lambda^{-1}(\cdot-p))|+|\nabla^2\beta(\lambda^{-1}(\cdot-p))|+
\lambda|\nabla\beta(\lambda^{-1}(\cdot-p))|+
\lambda^2,
\end{align}
and
\begin{align}
\label{5.71}
\Big|\epsilon\Big[A(\q),\epsilon\frac{\partial^2 A}{\partial
p_1^2}(\q)\Big]&\Big|\lesssim
\epsilon^2\lambda^{-2}|\nabla^2\beta(\lambda^{-1}(\cdot-p))|+\frac{\epsilon\lambda^2}{r^3(\lambda^2+r^2)}+
\epsilon|\nabla^2\beta(4\lambda^{-1}(\cdot-p))|\notag\\
&\quad+\epsilon\lambda|\nabla\beta(\lambda^{-1}(\cdot-p))|+\epsilon\lambda^2+
\frac{\epsilon}{r(\lambda^2+r^2)}|\nabla^2\beta(\lambda^{-1}(\cdot-p))|\notag\\
&\quad+\frac{\lambda^4}{r^4(\lambda^2+r^2)^2}+
\frac{\lambda^2}{r(\lambda^2+r^2)}|\nabla^2\beta(4\lambda^{-1}(\cdot-p))|\notag\\
&\quad+\frac{\lambda^3}{r(\lambda^2+r^2)}|\nabla\beta(4\lambda^{-1}(\cdot-p))|+\frac{\lambda^4}{r(\lambda^2+r^2)}.
\end{align}
From \eqref{5.69}--\eqref{5.71}, using $|a_{11}(\q)|\lesssim\epsilon^{3/2}$, one
 obtains
\begin{align}
\label{5.72}
a_{11}(\q)^4\int_{B^4\setminus B_{\lambda/4}(p)}\Big|\frac{\partial^2 A}{\partial p_1^2}(\q)\Big|^2\,dx&\lesssim\epsilon^6+\epsilon^4\int_{B^4\setminus B_{\lambda/4}(p)}\frac{\lambda^4}{r^6(\lambda^2+r^2)^2}\,dx+\epsilon^4\lambda^4\notag\\
&\quad+\epsilon^4\lambda^6+\epsilon^4\lambda^4\lesssim\epsilon^3,
\end{align}
\begin{align}
\label{5.73}
a_{11}(\q)^4\int_{B^4\setminus B_{\lambda/4}(p)}\Big|\nabla\Big(\frac{\partial^2 A}{\partial p_1^2}\Big)\Big|^2\,dx&\lesssim\epsilon^6\lambda^{-2}+\epsilon^6+\epsilon^4\int_{B^4\setminus B_{\lambda/4}(p)}\frac{\lambda^4}{r^8(\lambda^2+r^2)^2}\,dx+\epsilon^6\lambda^2\notag\\
&\quad+\epsilon^4\lambda^4\lesssim\epsilon^2,
\end{align}
and
\begin{align}
\label{5.74}
&\int_{B^4\setminus B_{\lambda/4}(p)}\Big|\epsilon\Big[A(\q),a_{11}(\q)^2\frac{\partial^2 A}{\partial p_1^2}(\q)\Big]\Big|^2\,dx\notag\\
&\lesssim\epsilon^8+\epsilon^6\int_{B^4\setminus B_{\lambda/4}(p)}\frac{\lambda^4}{r^6(\lambda^2+r^2)^2}\,dx+\epsilon^6\lambda^4+
\epsilon^6\lambda^6+\epsilon^6\lambda^4+\epsilon^6\int_{B_{2\lambda}(p)\setminus B_{\lambda}(p)}\frac{1}{r^2(\lambda^2+r^2)^2}\,dx\notag\\
&\quad+\epsilon^4\int_{B^4\setminus B_{\lambda/4}(p)}\frac{\lambda^8}{r^8(\lambda^2+r^2)^4}\,dx+\epsilon^4\int_{B_{\lambda/2}(p)\setminus B_{\lambda/4}(p)}\frac{\lambda^4}{r^2(\lambda^2+r^2)^2}\,dx\notag\\
&\quad+\epsilon^4\int_{B_{\lambda/2}(p)\setminus
B_{\lambda/4}(p)}\frac{\lambda^6}{r^2(\lambda^2+r^2)^2}\,dx+\epsilon^4\int_{B^4\setminus
B_{\lambda/4}(p)}\frac{\lambda^8}{r^2(\lambda^2+r^2)^2}\,dx\lesssim\epsilon^2.
\end{align}
Combining \eqref{5.72}--\eqref{5.74} yields
\begin{equation}
\label{5.75}
\Big\|a_{11}(\q)^2\frac{\partial^2 A}{\partial
p_1^2}(\q)\Big\|_{A(\q);1,2;B^4\setminus
B_{\lambda/4}(p)}\lesssim\epsilon\;.
\end{equation}
Since $\epsilon A(\q)=gI^1_{\lambda,p}g^{-1}$ and
$\epsilon\frac{\partial^2 A}{\partial p_1^2}=g\frac{\partial^2
I^1_{\lambda,p}}{\partial p_1^2}g^{-1}$ on $B_{\lambda/4}(p)$, one
has
\begin{align}
\Big|a_{11}(\q)^2\frac{\partial^2 A}{\partial
p_1^2}(\q)\Big|\lesssim\frac{\epsilon^2}{(\lambda^2+r^2)^{3/2}}\;,\notag\\
\Big|\epsilon\Big[A(\q),a_{11}(\q)^2\frac{\partial^2 A}{\partial
p_1^2}(\q)\Big]\Big|\lesssim\frac{\epsilon^2}{(\lambda^2+r^2)^2}\;,\notag\\
\Big|\nabla\Big(a_{11}(\q)^2\frac{\partial^2 A}{\partial
p_1^2}(\q)\Big)\Big|\lesssim\frac{\epsilon^2}{(\lambda^2+r^2)^2}\;.\notag\\
\notag
\end{align}
Hence,
\begin{align}
\int_{B_{\lambda/4}(p)}\Big|a_{11}(\q)^2\frac{\partial^2 A}{\partial
p_1^2}(\q)\Big|^2\,dx\lesssim\epsilon^4\int_{B_{\lambda/4}(p)}
\frac{1}{(\lambda^2+r^2)^3}\,dx\lesssim\epsilon^3,\label{5.76}\\
\int_{B_{\lambda/4}(p)}\Big|\nabla\Big(a_{11}(\q)^2\frac{\partial^2
A}{\partial
p_1^2}(\q)\Big)\Big|^2\,dx\lesssim\epsilon^4\int_{B_{\lambda/4}(p)}\frac{1}{(\lambda^2+r^2)^4}\,dx
\lesssim\epsilon^2,\label{5.77}\\
\int_{B_{\lambda/4}(p)}\Big|\epsilon\Big[A(\q),a_{11}(\q)^2\frac{\partial^2
A}{\partial
p_1^2}(\q)\Big]\Big|^2\,dx\lesssim\epsilon^2
\int_{B_{\lambda/4}(p)}\frac{1}{(\lambda^2+r^2)^4}\,dx\lesssim\epsilon^2.\label{5.78}
\end{align} From \eqref{5.76}-\eqref{5.78}, one obtains
\begin{equation}
\label{5.79}
\Big\|a_{11}(\q)^2\frac{\partial^2 A}{\partial
p_1^2}(\q)\Big\|_{A(\q);1,2;B_{\lambda/4}(p)}\lesssim\epsilon,
\end{equation}
and, finally, from \eqref{5.75} and \eqref{5.79}, the assertion of Lemma \ref{L3.10}
follows, for the case $i=1$, $j=1$. The remaining cases can be
proved similarly.\hfill$\Box$

\end{document}